\documentclass[letterpaper, 10 pt, conference]{ieeeconf}  

\usepackage{graphicx}

\usepackage{amsmath}
\usepackage{paralist}
\usepackage{amsfonts}
\usepackage{amssymb}
\usepackage{array}
\usepackage{amsthm}
\usepackage{subcaption}

\newtheorem{theorem}{Theorem}
\newtheorem{lemma}[theorem]{Lemma}
\newtheorem{corollary}[theorem]{Corollary}
\theoremstyle{definition}

\theoremstyle{remark}

\newcommand{\E}{\mathbb{E}}

\newcommand{\cvr}{\text{CVaR}_{\alpha}(c(x_{k},a))}
\newcommand{\Prb}{\mathbb{P}}

\usepackage{comment}
\usepackage{url}

\usepackage{cite}
\usepackage{accents}

\DeclareMathOperator*{\argmax}{argmax}

\title{\LARGE \bf
Controlled Information Fusion with Risk-Averse CVaR Social Sensors
}

\author{{\hspace{4cm} Sujay Bhatt and Vikram Krishnamurthy, {\em Fellow, IEEE}}
\thanks{S. Bhatt ({\tt\small sh2376@cornell.edu}) and V. Krishnamurthy ({\tt\small vikramk@cornell.edu}) are with the Department of Electrical and Computer Engineering, Cornell University, Ithaca, NY 14850. \newline
 This material is based upon work supported, in part by, the U. S. Army Research Laboratory and the U. S. Army Research Office under grant 12346080, National Science Foundation under grant 1714180, and Schmidt Sciences.}}

\usepackage{epstopdf}
\usepackage[font={small}]{caption}

\newsavebox{\smlmat}
\savebox{\smlmat}{$\begin{bmatrix}0.7&0.3\\0.3&0.7\end{bmatrix}$}
\IEEEoverridecommandlockouts

\begin{document}

\maketitle
\thispagestyle{empty}
\pagestyle{empty}

\begin{abstract}
Consider a multi-agent network comprised of risk-averse social sensors and a controller that jointly seek to estimate an unknown state of nature, given noisy measurements. 

The network of social sensors perform \textit{Bayesian social learning} - each sensor fuses the information revealed by previous social sensors along with its private valuation using Bayes' rule - to optimize a local cost function. The controller sequentially modifies the cost function of the sensors by discriminatory pricing (control inputs) to realize long term global objectives. 

We formulate the stochastic control problem faced by the controller as a Partially Observed Markov Decision Process (POMDP) and derive structural results for the optimal control policy as a function of the risk-aversion factor in the Conditional Value-at-Risk (CVaR) cost function of the sensors. We show that the optimal price sequence when the sensors are risk-averse is a \textit{super-martingale}; i.e, it decreases on average over time.  

\end{abstract}

\section{INTRODUCTION}

A social sensor is an information processing system having the following attributes:
\begin{compactenum}
\item[i.)] It affects the behaviour of other sensors.
\item[ii.)] It shares quantized information (decisions/actions) and has its own dynamics.
\item[iii.)] It has limited processing capabilities - boundedness.
\item[iv.)] It is rational - fuses all available information using Bayes' rule to take optimal action.
\end{compactenum}
Social learning is the process by which social sensors are influenced by the behaviour of other sensors in a multi-agent network. 
We present a model of Bayesian social learning with focus on understanding the interaction between a controller and a multi-agent network of risk-averse social sensors. 

\cite{Kri12} considers the interaction of a controller and a network of social sensors, where the sensors perform social learning to estimate an unknown parameter and optimize a local utility function. The controller seeks to detect a change in the parameter as soon as possible by observing the actions of the sensors. This is the well studied \textit{Controlled Sensing} problem, see \cite{Kri16}, where the observation (action) statistics is controlled to meet the desired objectives. In \cite{Kri12}, well known inefficiencies of the standard social learning model \cite{BHW92,Ban92} like herding (sensors choose the same action irrespective of their private information) and informational cascades (information fusion results in no improvement in uncertainty) are shown to be the consequence of the belief state (probability distribution on the parameter) belonging to suitably defined regions in the belief space. Using this characterization, structural results (threshold policies) on the optimal policy of the controller are derived in \cite{Kri12}. \cite{KB16} extends the analysis to the case where the sensors display an aversion to risk; i.e, social sensors having risk aversion as an additional attribute. It is shown that when the sensors are risk-averse (modeled using CVaR), the herding behaviour is more pronounced - social learning is absent when the sensors are sufficiently risk-averse. In this paper, we extend the results in \cite{KB16} to the case where the controller directly modifies the cost function of the social sensors, and this problem is termed as the \textit{controlled information fusion} problem. The controller fuses (aggregates) the information on the state revealed in the form of decisions of the risk-averse social sensors.

The contribution of this paper is two fold:
\begin{compactenum}
\item  We study the interaction of a controller and a network of risk-averse social sensors, using a Partially Observed Markov Decision Process (POMDP) framework. Unlike controlled sensing, the interaction is such that the controller can influence the information fusion in social sensors.
\item We obtain structural results for the optimal policy of the controller and characterize the properties of the optimal (price) sequence, when the social sensors are risk-averse.
\end{compactenum}
\cite{BOOV06,BOOV08} consider monopoly pricing in the presence of social learning and establish various properties of the value function and the optimal policy for the monopoly. It is shown that using discriminatory pricing the monopoly is able to delay the process of herding in risk-neutral social sensors, to suit its needs. The optimal price (control) sequence is shown to be a super-martingale. We consider monopoly pricing and social learning under CVaR risk-measure{\footnote{A risk measure $\varrho : \mathcal{L} \rightarrow \mathbb{R}$ is a mapping from the space of measurable functions to the real line which satisfies the following properties: (i) $\varrho(0)=0$. (ii) If $S_{1}, S_{2} \in \mathcal{L}$ and $S_{1} \leq S_{2} ~\text{a.s}$ then $\varrho(S_{1}) \leq \varrho(S_{2})$. (iii) if $a\in\mathbb{R}$ and $S\in\mathcal{L}$, then $\varrho(S+a) = \varrho(S)+a $. The risk measure is coherent if in addition $\varrho$ satisfies: (iv) If $S_{1}, S_{2} \in \mathcal{L}$, then $\varrho(S_{1}+S_{2}) \leq \varrho(S_{1}) + \varrho(S_{2})$. (v) If $a \geq 0$ and $S\in\mathcal{L}$, then $\varrho(aS)=a\varrho(S)$. The expectation operator is a special case where subadditivity is replaced by additivity.}}; see \cite{MJ10} for an overview of risk measures. 

CVaR is an extension of VaR that gives the total loss given a loss event, and is a coherent risk measure; see \cite{RU00}. The value at risk (VaR) is the percentile loss namely,  $\text{VaR}_\alpha(x) = \min\{z: F_x(z) \geq \alpha\} $ for cdf $F_x$, and $\text{CVaR}_\alpha(x) = \E\{X | X > \text{VaR}_\alpha(x)\}$.
CVaR is one of the `big' developments in risk modelling in the last 15 years. 
In this paper, we choose CVaR risk measure as it exhibits the following properties:
(i) It associates higher risk with higher cost.
(ii) It ensures that risk arises only from the services.
(iii) It is convex. 
\subsection*{Organization and Main Results}
Sec.~\ref{sec:PPCO} details the Bayesian social learning model for the process of information fusion by CVaR social sensors and the pricing protocol employed by the controller. The controller's long term objective is to maximize a discounted reward function.\\
Sec.~\ref{sec:SOPP} formulates the stochastic control problem faced by the controller as a Partially Observed Markov Decision Process (POMDP) and is solved using dynamic programming. The structure of the value function and optimal policy is completely characterized.\\
Sec.~\ref{sec:POPS} describes the nature of the price sequence that is input to the multi-agent network. It is shown that the controller prices high initially and subsequently lowers the price, i.e the price input sequence over time is a super-martingale.

\section{CVaR Social Learning Model and Controller Objective} \label{sec:PPCO}
We consider the classical sequential social learning framework \cite{Cha04,Kri12,KB16}. The social sensors and the controller jointly seek to estimate an unknown state of nature to meet the desired objectives. The controller sequentially chooses price inputs to the multi-agent network in exchange for services and the sensors decide to utilize the services depending on its quality. The decision (action) of each sensor depends on the cost, risk factor, private valuation and the decisions of the other sensors in the network{\footnote{Amazon Web Services (AWS), for example, provides an on-demand cloud platform with a wide range of services like storage, developer tools, analytics etc for client-side applications at different prices. AWS is only partially aware of the \textit{quality} of its services, while the clients learn about it from the experience of other clients and self-valuation.}}. 

Let $x  \in \mathcal{X} = \lbrace 1(\text{Low}),2(\text{High})\rbrace$ denote the state. In this paper, we study the problem of localization; i.e, the quality is a random variable. Let the initial distribution (on the quality) be denoted as $\pi_{0} = (\pi_{0}(i),i\in \mathcal{X})$, where $\pi_{0}(i) = \Prb(x_{0}=i)$.

Each sensor acts once in a predetermined sequential order indexed by $k = 1,2,\cdots$. The index $k$ can also be viewed as the discrete time instant when sensor $k$ acts. Each sensor $k$ obtains noisy private valuations, $y_k \in \mathcal{Y} = \{1(\text{Low}),2(\text{High})\}$, of the quality $x_k$ and considers this in addition to the actions of its predecessors. The controller does not have any information about $x_k \in \mathcal{X}$ but infers it from the information revealed by the actions of the individual sensors, $a_k \in \mathcal{A} = \lbrace 1(\text{Don't Utilize}),2(\text{Utilize}) \rbrace$, and chooses the price inputs{\footnote{The range of prices chosen by the controller is normalized to $[0,1]$ for convenience.}}  $u_k \in [0,1]$ at each time $k$ (or at each sensor $k$).  

\subsection{CVaR Social Learning Model and Pricing Protocol}

\begin{figure}[!t] 
\centering
{\hspace{-0cm}}\includegraphics[scale=0.35]{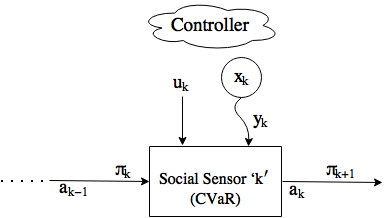}
\caption{CVaR Social Learning model. The social sensor $k$ receives the public belief $\pi_{k-1}$ from all its predecessors. $y_k$ denotes the private valuation of the quality and $u_k$ denotes the price charged by the controller for the services. The decision $a_k$ is shared by the controller and the updated public belief $\pi_{k+1}$ is received by the successive sensors.}
\label{SL}
\end{figure}

The social learning model and the pricing protocol of the controller is as follows:
\begin{compactenum}
\item[1.] \textit{Social Sensor's Private Observation}: Social sensor $k$'s private observation denoted by $y_{k} \in \mathcal{Y} = \{1,2\}$ is a noisy measurement of the true quality. It is obtained from the observation likelihood distribution as,
\begin{equation}\label{eq:obs_m}
B_{ij} = \Prb(y_{k}=j|x_{k}=i)
\end{equation} 
The discreteness of the observation distribution captures the \textit{boundedness} or the limited processing capabilities of the sensor.
\item[2.] \textit{Social Learning and Private Belief update}: Social sensor $k$ updates its private belief by fusion of the observation $y_{k}$ and the prior public belief $\pi_{k-1}(i) = \Prb(x_k=i|a_{1},\hdots,a_{k-1})$ as the following Hidden Markov Model (HMM) update
\begin{equation} \label{eq:PBU}
\eta^{y_k}_{k} = \frac{B_{y_{k}}\pi_{k-1}}{\textbf{1}'B_{y_{k}}\pi_{k-1}}
\end{equation}
where $B_{y_k}$ denotes the diagonal matrix having $[\Prb(y_{k}|x_{k}=~1) ~ \Prb(y_{k}|x_{k}=2))]$ along the diagonal and $\textbf{1}$ denotes the $2$-dimensional vector of ones. HMM update is a consequence of Bayes' rule, information on the state conditioned on the new observation.
\item[3.] \textit{Social Sensor's Action}: Social sensor $k$ executes an action $a_{k}\in\mathcal{A}=\lbrace 1,2\rbrace$ to myopically minimize its cost. Let $c(x_k,a_{k})$ denote the cost incurred if the sensor takes action $a_{k}$ when the underlying state is $x_k$. 

The form of the state-action dependent cost is taken as $c(x_k,a_{k}) = u_k - v(x_k)$ (see \cite{Cha04} for a justification), where $v$ is the valuation of the services by each sensor and $u_k$ is the price chosen by the controller at time $k$. It is assumed without loss of generality that 
\[ v(x_k) = \left\{ \begin{array}{ll}
         0 & \mbox{if $x_k = 1$};\\
        1 & \mbox{if $x_k = 2$}.\end{array} \right. \] 
The state-action dependent costs for $x  \in \mathcal{X}$ are thus given as:
\[ c(x_k,a_{k}) = \left\{ \begin{array}{ll}
         \begin{bmatrix}
0 \\ 0
\end{bmatrix} & \mbox{if $a_k = 1$};\\
       \begin{bmatrix}
u_k \\ u_k-1
\end{bmatrix} & \mbox{if $a_k = 2$}.\end{array} \right. \]         

The sensor chooses an action $a_{k}$ to minimize the CVaR measure as
\begin{align}
a_{k} &=  {\underset{a \in \mathcal{A}}{\text{argmin}}} \{ \cvr \} \\
&=  {\underset{a \in \mathcal{A}}{\text{argmin}}} \{ {\underset{z \in \mathbb{R}}{\text{min}}} ~ \{ z + \frac{1}{\alpha} \mathbb{E}_{y_{k}}[{\max} \{ (c(x_{k},a)-z),0 \rbrace] \} \} \nonumber
\end{align}
Here $\alpha \in (0,1]$ reflects the degree of risk-aversion for the sensor (the smaller $\alpha$ is, the more risk-averse the sensor is). Note that when $\alpha = 1$, the cost function is the risk-neutral cost function as in \cite{Cha04,BOOV06,BOOV08}.
Define 
\begin{equation}
\hspace{-0.5cm} \mathcal{G}_{k} := \sigma \text{- algebra generated by}~ (u_1,a_{1},u_2,a_{2},\hdots,u_k,y_{k}) 
\end{equation} 
$\mathbb{E}_{y_{k}}$ denotes the expectation with respect to private belief, i.e, $\mathbb{E}_{y_{k}} = \mathbb{E}[.|\mathcal{G}_{k}]$ when the private belief is updated after observation $y_{k}$.
 
\item[4.] \textit{Controller Reward}: We consider two possible reward functions for the controller. The controller chooses one of the following reward functions at $k=0$ and accrues the corresponding reward at each time $k$ as
\begin{compactenum}
\item[Case 1.] \underline{Self-Interested}: The controller accrues a reward when the sensors utilize the services,
\begin{equation} \label{eq:PM}
r_{u_k} = (u_{k}-\beta) \mathcal{I}(a_k=2| \pi_k). 
\end{equation}
\item[Case 2.] \underline{Altruistic}: The controller accrues a reward when the sensors act according to their valuations,
\begin{equation} \label{eq:AM}
r_{u_k} = (u_{k}-\beta) \mathcal{I}(a_k=y_k | \pi_k). 
\end{equation}
\end{compactenum} 
Here $\mathcal{I}$ denotes the indicator function and $\beta \in (0,1)$ is a fixed{\footnote{Note that $\beta$ could be made state dependent without affecting the nature of the results in the paper. Here it is assumed to be independent of the state for simplicity.}} cost incurred by the controller. It could denote the cost of service. The controller being self-interested can be seen as profit maximizing, and being altruistic can be seen as social welfare maximizing{\footnote{Social welfare is maximized when the controller and the sensors in the multi-agent network take decisions considering network externalities; see (Chapter~$4$,~\cite{Cha04}). We shall see in Sec.~\ref{sec:POPS} that $ \mathcal{I}(a=y)$ improves the value of information fused by the successive sensors, thereby promoting welfare.}}.

\item[5.] \textit{Public Belief update}: Sensor $k$'s action is shared by the controller with the multi-agent network  and the public belief on the quality is updated according to the social learning Bayesian filter (see \cite{Kri12,KB16}) as follows 
\begin{equation} \label{eq:SLF}
\pi_{k} = T^{\pi} (\pi_{k-1},a_k)  = \frac{R_{a_{k}}^{\pi_{k-1}}\pi_{k-1}}{\textbf{1}'R_{a_{k}}^{\pi_{k-1}}\pi_{k-1}}.
\end{equation}
Here, $\pi_{k}(i) = \Prb(x_k=i|a_{1},\hdots,a_{k})$, $R_{a_{k}}^{\pi_{k-1}} = \text{diag}(\Prb(a_{k}|x=i,\pi_{k-1}),i \in \mathcal{X})$, where $\Prb(a_{k}|x=i,\pi_{k-1}) = {\underset{y \in \mathcal{Y}}{\sum}} \Prb(a_{k}|y,\pi_{k-1})\Prb(y|x_{k}=i)$ and 
\begin{equation*} 
{\hspace{-0.5cm}} \Prb(a_{k}|y,\pi_{k-1}) = \left\{ \begin{array}{ll}
         1 & \mbox{if $ a_{k} = {\underset{a \in \mathcal{A}}{\text{argmin}}} \{ \cvr \}$} ; \\
         0 & \mbox{$\text{otherwise}$}.\end{array} \right. 
\end{equation*}          
Note that $\pi_k$ belongs to the unit simplex  
{\small{
\begin{equation*}
{\hspace{-0.5cm}} \Pi(2){\overset{\Delta}{=}}\lbrace \pi \in \mathbb{R}^{2} : \pi(1)+\pi(2) = 1, 0 \leq \pi(i) \leq 1 ~\text{for} ~ i \in \{1,2\} \rbrace 
\end{equation*}}}
 Here the expectation in CVaR measure is with respect to the sigma-algebra $\mathcal{G}_{k}$. Social learning filter update is a consequence of Bayes' rule, information on the state conditioned on the new action.
\item[6.]\textit{Information Fusion Price}: Let the history recorded by the controller and the multi-agent network be denoted as $\mathcal{H}_k = \{\pi_0,u_1,a_1,\cdots,u_k,a_k \}$. The controller chooses $u_{k+1} =~\mu_{k+1}(\mathcal{H}_k) \in [0,1]$ for the sensor $k+1$ and the protocol is repeated for all the sensors in the system. Here $\mu_{k+1}$ denotes the pricing policy at time $k+1$. 
\end{compactenum}
Fig.~\ref{SL} shows the CVaR social learning model.

\subsection{Controlled Fusion Objective}
The controller chooses the price inputs to the social sensors sequentially as
\begin{equation} \label{eq:con_H}
u_{k} =~\mu_{k}(\mathcal{H}_{k-1}) \in [0,1]
\end{equation}
where $\mathcal{H}_k = \{\pi_0,u_1,a_1,\cdots,u_{k-1},a_{k-1} \}$. Since $\mathcal{H}_k$ is increasing with time $k$, to implement a controller, it is useful to
obtain a sufficient statistic that does not grow in dimension. The public belief $\pi_{k-1}$ computed via the social learning filter (\ref{eq:SLF}) forms a sufficient statistic for $\mathcal{H}_k$ and{\footnote{The rewards are a function of the price and the state (see Lemma~\ref{lem:RM} and Lemma~\ref{lem:WM}), and hence restriction to Markov policies is without loss of generality.}} (\ref{eq:con_H}) can be written as
\begin{equation}
u_{k} = \mu_k(\pi_{k-1}).
\end{equation}
The controller maximizes the cumulative discounted reward
\begin{equation} \label{eq:RM}
J_{\mu}(\pi) = \E_{\mu} \{ \sum_{k=1}^{\infty} \rho^k r_{u_{k}}| \pi_0 = \pi \}.
\end{equation}
Here $u_k = \mu(\pi_{k-1})$ and $\rho \in [0,1)$ denotes the economic discount factor indicating the degree of impatience of the controller. 
In (\ref{eq:RM}), the controller seeks to find the optimal stationary policy $\mu^*$ such that
\begin{equation} \label{eq:Opt_J}
J_{\mu^*}(\pi_0) = {\text{sup}}_{\mu \in \boldsymbol{\mu}} J_{\mu}(\pi_0).
\end{equation}
The stochastic control problem faced by the controller is formulated as a partially observed Markov decision process (POMDP) with dynamics (\ref{eq:SLF}) and objective (\ref{eq:Opt_J}), and is solved using dynamic programming. 

\section{Structure of Optimal Pricing Policy} \label{sec:SOPP}
In this section, we characterize the nature of the optimal pricing policy for (\ref{eq:RM}). It is shown that due to the structure of the social learning filter in (\ref{eq:SLF}), the choice of price inputs reduces from a continuum to a finite number at every belief.
\subsubsection*{Assumptions}
\begin{compactenum}
\item[(A1)] The observation distribution $B_{xy} = \mathbb{P}(y|x)$ is TP2(total positive of order 2), i.e, the determinant of the matrix $B$ is non-negative; see \cite{Kri16}.
\end{compactenum}
The optimal policy $\mu^*$ and the value function $V(\pi)$ for the POMDP satisfy the Bellman's dynamic programming equation
\begin{align} \label{eq:RMVP}
Q(\pi,u) &= r_u + \rho \sum_{a \in \mathcal{A}} V(T^{\pi} (\pi,a)) \sigma(\pi,a), \nonumber \\  
\mu^*(\pi) &= \argmax_{u \in [0,1]} Q(\pi,u),  \nonumber \\ 
V(\pi) &= \max_{u \in [0,1]} Q(\pi,u), ~~ J_{\mu^*}(\pi_0) = V(\pi_0).
\end{align}

\begin{theorem} \label{thm:RM}
Given a risk-aversion factor $\alpha \in (0,1]$, let $u^{H}(\pi) = 1-\frac{\eta^{y=2}(1)}{\alpha}$ and $u^{L}(\pi) = 1-\frac{\eta^{y=1}(1)}{\alpha}$ denote two possible prices at the belief $\pi$. Under (A1), the $Q$ function (\ref{eq:RMVP}) can be simplified for the rewards (\ref{eq:PM}) and (\ref{eq:AM}) as
\begin{compactenum}
\item[Case 1.)] Self-Interested:\\
\begin{equation*} 
\hspace{-0.8cm} Q(\pi,u) = \left\{ \begin{array}{lll}
         (u-\beta) + \rho V( \pi) & \mbox{if $u \in [0,u^{L}(\pi)]$};\\
      (u-\beta) \times \textbf{1}^\prime B_{y=2} \pi  \\ + \rho \E V(\pi) & \mbox{if $u \in (u^{L}(\pi),u^{H}(\pi)]$}; \\
       0 & \mbox{$\text{otherwise}$}.\end{array} \right. 
 \end{equation*}  
and $V(\pi) = \max Q(\pi,u)$, where $V(\pi) \geq 0$.  
\item[Case 2.)] Altruistic:\\
\begin{equation*}
{\hspace{-0.8cm}}Q(\pi,u) = \left\{ \begin{array}{ll}
              (u-\beta) + \rho \E V(\pi) & \mbox{if $u \in (u^{L}(\pi),u^{H}(\pi)]$}; \\
       0 & \mbox{$\text{otherwise}$}.\end{array} \right. 
        \end{equation*} 
and $V(\pi) = \max Q(\pi,u)$, where $V(\pi) \geq 0$.       
\end{compactenum}         
\end{theorem}
Here, $$\E V(\pi) = \textbf{1}'B_{y=1}^{\pi} \pi \times V(\eta^{y=1}) + \textbf{1}'B_{y=2}^{\pi} \pi \times V(\eta^{y=2}).$$ \\ 
The prices $u^{H}(\pi)$ and $u^{L}(\pi)$ are such that sensors utilize the services when $y=2$ and $y=\{1,2\}$ respectively. The proof of Theorem~\ref{thm:RM} will be given in the appendix.  Theorem~\ref{thm:RM} represents the Q function (\ref{eq:RMVP}) over a price input range $[0,1]$ for rewards (\ref{eq:PM}) and (\ref{eq:AM}) respectively, in \textit{three} and \textit{two} regions. The following corollaries highlight why such partitions are useful. 
\begin{corollary} \label{cor:FnP}
Let the controller reward be given by (\ref{eq:PM}). At every public belief $\pi \in \Pi(2)$, it is sufficient to choose one of the three prices $\{ u^{L}(\pi), u^{H}(\pi),u^{H}(\pi)+\epsilon\}$ for any $\epsilon>~0$. \qed
\end{corollary} 
\begin{corollary} \label{cor:FnA}
Let the controller reward be given by (\ref{eq:AM}). At every public belief $\pi \in \Pi(2)$, it is sufficient to choose one of the two prices $\{u^{H}(\pi),u^{H}(\pi)+\epsilon\}$ for any $\epsilon>~0$. \qed
\end{corollary}
The following theorem completely characterizes the optimal pricing policy when the controller aims to maximize the reward. The proof is given in the appendix.
\begin{theorem} \label{thm:OPRM}
For every public belief $\pi \in \Pi(2)$ and an $\epsilon>~0$, the optimal policy $\mu^*(\pi) = \argmax_u Q(\pi,u)$ is given as
\begin{compactenum}
\item[Case 1.)] Self-Interested:
\begin{equation} \label{eq:OPRM}
\mu^*(\pi) = \left\{ \begin{array}{lll}
         u^{H}(\pi) + \epsilon & \mbox{if $\pi(2) \in [0,\pi^*(2))$};\\
       u^{H}(\pi) & \mbox{if $\pi(2) \in [\pi^*(2), \pi^{**}(2))$}; \\
       u^{L}(\pi)  & \mbox{$\pi(2) \in [\pi^{**}(2), 1]$}.\end{array} \right. 
        \end{equation}   
   for $\pi^*(2),~\pi^{**}(2) \in [0,1]$.  
\item[Case 2.)] Altruistic:
\begin{equation} \label{eq:OPWM}
\mu^*(\pi) = \left\{ \begin{array}{lll}
         u^{H}(\pi) + \epsilon & \mbox{if $\pi(2) \in [0,\hat{\pi}^*(2))$};\\
       u^{H}(\pi) & \mbox{if $\pi(2) \in [\hat{\pi}^*(2), 1]$}.\end{array} \right. 
        \end{equation}   
   for  $\hat{\pi}^*(2) \in (0,1)$.  
   \end{compactenum}
\end{theorem}
\subsubsection*{\underline{Discussion}}
\begin{figure}
\centering
\begin{subfigure}{.25\textwidth}
  {\hspace{-0.52cm}}\includegraphics[scale=0.25]{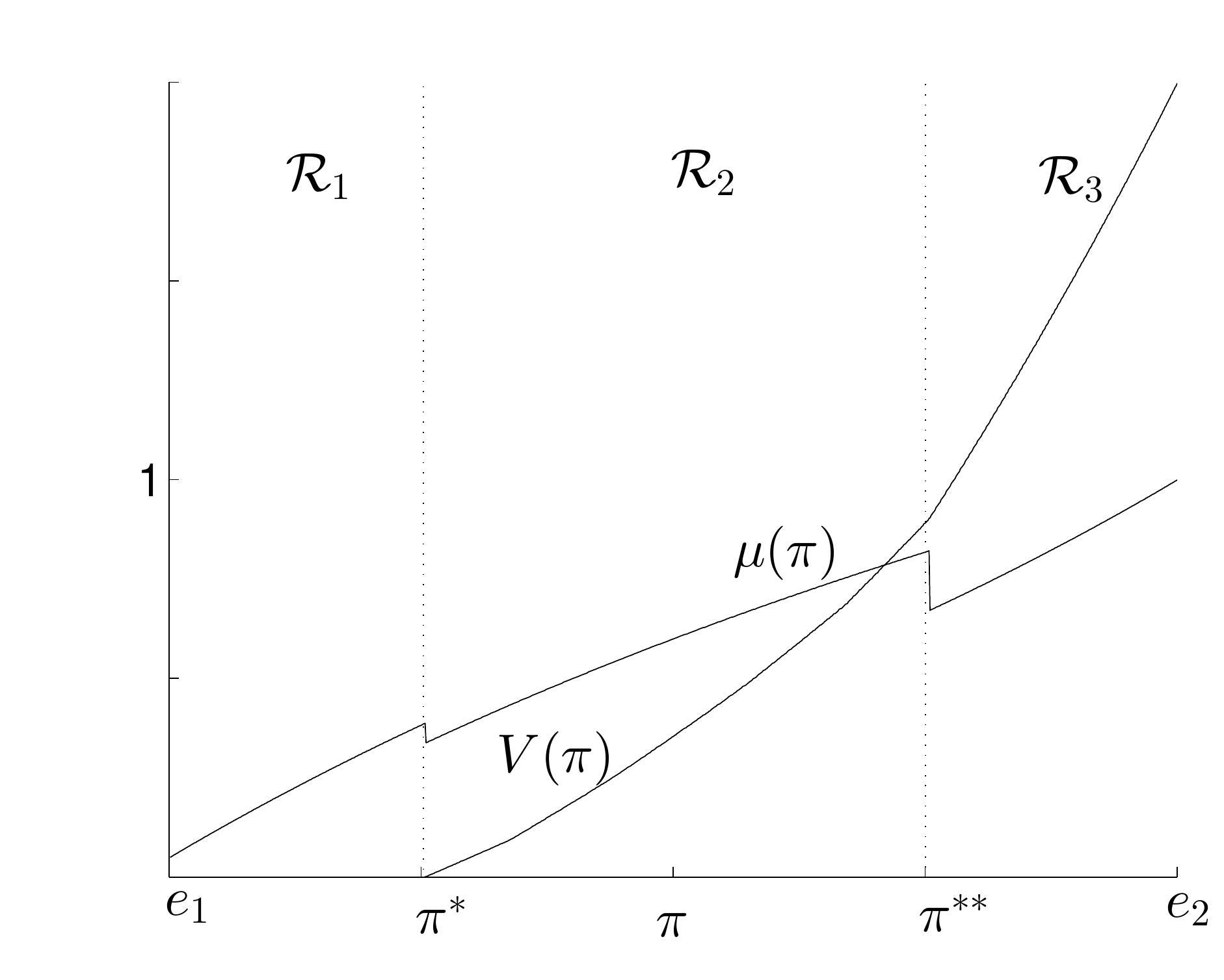}
  \caption{Risk-aversion factor $\alpha=0.9$ }
\label{RM1}
\end{subfigure}%
~
\begin{subfigure}{.25\textwidth}
 {\hspace{-0.75cm}} \includegraphics[scale=0.25]{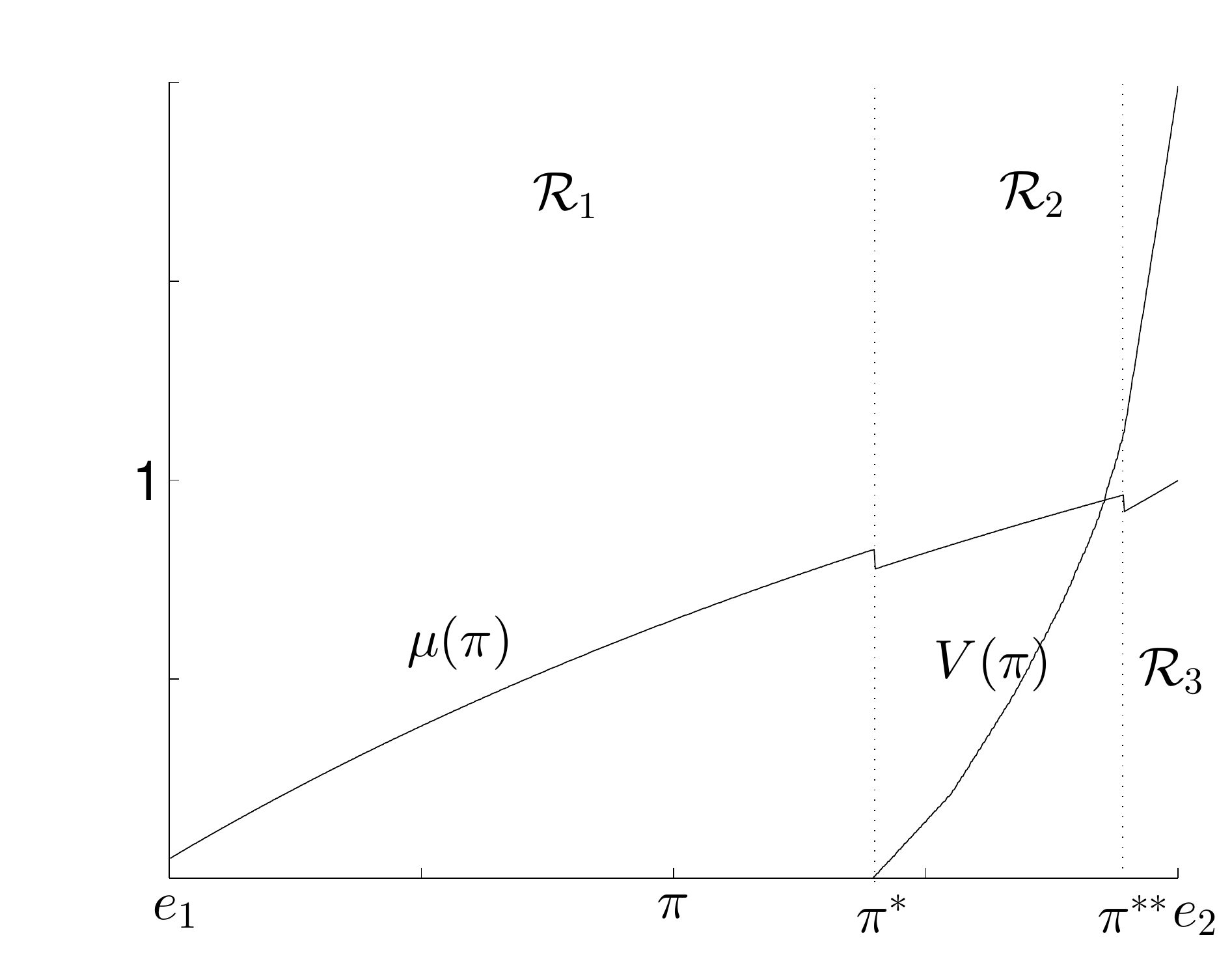}
 \caption{$\alpha=0.3$ }
\label{RM2}
\end{subfigure}
\caption{Value function and optimal pricing policy of the controller in the \textit{self-interested} case. $B=$~\usebox{\smlmat}, the discount factor $\rho = 0.7$, and $\mathcal{R}_1 - [0, \pi^*),~\mathcal{R}_2 - [\pi^*, \pi^{**}),~\text{and}~\mathcal{R}_3 - [\pi^{**}, 1]$ are the cut-off, social learning and herding regions respectively. It can be seen that the width of $\mathcal{R}_1$ increases with increased aversion to risk. This is equivalent to saying that risk-averse sensors that show an increased aversion to risk, choose to utilize the services only when they are reasonably certain about the quality. So it is profitable to the controller if it offers services only when it believes that the quality is high.}
\end{figure}

\
From Theorem~\ref{thm:RM}, Corollary~\ref{cor:FnP} and Corollary~\ref{cor:FnA}, the value function in (\ref{eq:RMVP}) can be represented as 
\begin{align*}
 (\text{Self-Interested}) \\
V(\pi) &= \max \{0, (u^{L}(\pi)-\beta) + \rho V(\pi), \\
  &(u^{H}(\pi)-\beta) \times \textbf{1}^\prime B_{y=2} \pi   + \rho \E V(\pi) \} \\
{\hspace{-1cm}} (\text{Altruistic}) \\
V(\pi) &= \max \{0, (u^{H}(\pi)-\beta) +  \rho \E V(\pi) \}
\end{align*}
The key takeaway is that due to the structure of the social learning filter, the choice of price inputs at every belief is reduced to a finite number of values instead of the range $[0,1]$. Characterizing the optimal policy amounts to selecting among the these price inputs as a function of the public belief. Theorem~\ref{thm:OPRM} completely determines the regions in the belief space $\Pi(2)$ where it is optimal to choose a particular price input. Fig.~\ref{RM1} and Fig.~\ref{RM2} show the value function and the optimal policy for two different risk-aversion factors ($\alpha$) in a simple numerical example.

Let $\mathcal{R}_1,~\mathcal{R}_2,~\mathcal{R}_3$ denote the three regions determined by Theorem~\ref{thm:OPRM} where $u^{H}(\pi)+\epsilon,~u^{H}(\pi),~u^{L}(\pi)$ respectively are optimal. $\mathcal{R}_1$ is the \textit{cut-off} region - the controller terminates the services to the multi-agent network. In the \textit{Self-Interested} case, the price inputs are such that no sensor has an incentive to utilize the services. $\mathcal{R}_2$ is the \textit{social learning} region - the sensors act according to their private valuations. The price inputs are such that the sensor having a high valuation $y=2$ will utilize the services, while the sensor having low valuation $y=1$ finds it prohibitive. Since the sensors act according to their valuation, sensor deciding at a future instant can successfully infer the private valuation of its predecessors; in other words, the information fusion reduces uncertainty about the quality of the service. $\mathcal{R}_3$ is the \textit{herding} region - every sensor utilizes the services. The controller chooses a low input $u^{L}(\pi)(< u^H(\pi))$, which prompts the sensor with even a low valuation $y=1$ to utilize the service. Notice that when the controller chooses $u = u^L(\pi)$ (when $\pi(2) \in [\pi^{**}(2),1]$), the value function is $V(\pi) = \frac{(u^L(\pi)-\beta)}{(1-\alpha)}$ - a fixed payoff. This means the controller induces a herd (sensors choose the same action irrespective of their private valuation) that leads to an information cascade (information fusion results in no improvement in uncertainty) - public belief is frozen. 

In the \textit{Altruistic} case, the price inputs (two at every belief) are chosen so as to encourage the sensors to act according to their valuations. This implies that the controller chooses inputs to maximize the width of the \textit{social learning} region $\mathcal{R}_2$. The \textit{herding} region $\mathcal{R}_3$ is absent as $u = u^L(\pi)$ is not chosen by the controller. The \textit{cut-off} region indicates the flexibility to terminate the services when the expected valuation is less than the cost of service.
\section{Properties of the Optimal Price Sequence} \label{sec:POPS}
In this section, we describe the relation between the optimal policy (\ref{eq:OPRM}) and (\ref{eq:OPWM}), and the price sequence $u_k = \mu^*(\pi_k)$. 

\begin{theorem} \label{thm:MPT}
Let $\mathcal{F}_k$ be the $\sigma$-algebra generated by $(u_1,a_{1},u_2,a_{2},\hdots,u_{k-1},a_{k-1},u_k,a_k)$, where $\pi_0$ is the initial belief. The optimal price sequence~$u_k = \mu^*(\pi_{k-1})$ is a super-martingale{\footnote{Decreases on average over time.}} when the quality is a random variable for any $\alpha \in (0,1]$. 
\end{theorem}
\subsubsection*{\underline{Discussion}}\
When the controller is profit maximizing or self-interested, it initially chooses higher price inputs to encourage sensors with higher valuation to utilize the services. Decisions at higher prices are more informative{\footnote{Informativeness is in the sense of Blackwell; see \cite{Kri16}. For any two observation matrices $B_1$ and $B_2$, $B_1$ is more informative than $B_2$ in the Blackwell sense ($B_1 \succ_B B_2$) if $B_2 = B_1Q$, for any stochastic matrix $Q$. Note here that when $u=u^H(\pi)$, the action likelihood matrix in (\ref{eq:SLF}) $R^H=B$; and when $u=u^L(\pi)$, the action likelihood matrix $R^L=\begin{bmatrix}
1 & 0\\
1 & 0
\end{bmatrix}$. We have for $Q = \begin{bmatrix}
1 & 0\\
1 & 0
\end{bmatrix},~$ $R^L = R^H Q \Rightarrow R^H \succ_B R^L$.}}, which in turn results in higher public beliefs when $a=2$. Due to the concavity of the pricing policy, higher belief causes the future price inputs to increase. 
Once sufficient information about the quality is accumulated, the controller either chooses low price inputs to allow every sensor to utilize the services or terminates its services to the multi-agent network. 

When the controller is altruistic, it always chooses high price inputs to encourage the sensors to act according to their private valuations.

\appendix
\begin{lemma}[\cite{Kri16}] \label{lem:IncPB}
Let $\eta^y$ denote the private belief update~(\ref{eq:PBU}) with a public prior belief $\pi$. Under (A1), $\eta^y$ is increasing{\footnote{$\pi_{2}\geq\pi_{1}$ if the determinant
\[ 
\begin{vmatrix} 
\pi_1(1) & \pi_1(2) \\ 
\pi_2(1) & \pi_2(2) 
\end{vmatrix}  \geq 0 \]}} in~$y$, i.e, $\eta^{y=1}(1) >\eta^{y=2}(1)$. \qed
\end{lemma}
\begin{theorem}[\cite{Kri16}] \label{thm:VFRM}
Let the instantaneous rewards be non-decreasing in $\pi$. Under (A1), the value function $V(\pi)$ with finite number of actions at every belief, is monotone and convex. \qed
\end{theorem}
 
\begin{lemma} \label{lem:RM}
The instantaneous reward $(u-\beta) \mathcal{I}(a=2|\pi)$ is given as
\begin{equation} \label{eq:InsRW}
\sum_{j \in \mathcal{Y}} \sum_{i \in \mathcal{X}} (u-\beta) \mathcal{I}(u \leq 1-\frac{\eta^{y=j}(1)}{\alpha}) B_{ij} \pi(i). 
\end{equation} 
\end{lemma}
\begin{lemma} \label{lem:WM}
The instantaneous reward $(u-\beta) \mathcal{I}(a=y|\pi)$ is given as
\begin{equation}
 (u-\beta)~\mathcal{I}(u^{L}(\pi) < u \leq u^{H}(\pi)). 
\end{equation} 
\end{lemma}
The proofs follow from the structure of the social learning filter (see Theorem~$2$, \cite{KB16}), property of the CVaR measure (see Lemma~$6$, \cite{KB16}), and Bayes' rule. It is omitted. \\
We will prove Theorem~\ref{thm:RM} and Theorem~\ref{thm:OPRM} for the \textit{Self-Interested} case. The proof for the \textit{Altruistic} case follows similarly. \\
\textit{\underline{Proof of Theorem~\ref{thm:RM}}}: \\  Consider $Q(\pi,u)$ as in (\ref{eq:RM}) for $u \in [0,1]$. 
\begin{compactenum}
\item[i.)] Let $u \in [0, u^{L}(\pi)]$. Recall that $u^{L}(\pi) = 1 - \frac{\eta^{y=1}(1)}{\alpha}$. The instantaneous reward in (\ref{eq:InsRW}) is~$(u-\beta)$. The continuation payoff $\sum_{a \in \mathcal{A}} V(T^{\pi} (\pi,a)) \sigma(\pi,a)$ is given as follows. From (\ref{eq:SLF}), $R_{a}^{\pi}    =  \begin{bmatrix}
0 & 1 \nonumber \\
0 & 1 \end{bmatrix}$.
\begin{align} \label{eq:P1}
(\Rightarrow) \sum_{a \in \mathcal{A}} V(T^{\pi} (\pi,a)) \sigma(\pi,a) &= V(\pi). \nonumber \\
\therefore Q(\pi,u) = (u-\beta) &+ \rho V(\pi).
\end{align}
\item[ii.)] Let $u \in (u^{L}(\pi), u^{H}(\pi)]$. The instantaneous reward in (\ref{eq:InsRW}) is~$(u-\beta)\times \textbf{1}^\prime B_{y=2} \pi $. From (\ref{eq:SLF}), $R_{a}^{\pi}  = B$.
\begin{align} \label{eq:P2}
(\Rightarrow) \sum_{a \in \mathcal{A}} V(T^{\pi} (\pi,a)) \sigma(\pi,a) &= \E V(\pi).  \nonumber \\
\therefore Q(\pi,u) = (u-\beta) \times \textbf{1}^\prime B_{y=2} \pi   &+ \rho \E V(\pi).
\end{align}
 \item[iii.)] Let $u \in (u^{H}(\pi), 1]$. This implies that $u > 1- \frac{\eta^{y=2}(1)}{\alpha}$. The instantaneous reward in (\ref{eq:InsRW}) is $0$. From (\ref{eq:SLF}), $R_{a}^{\pi}    =  \begin{bmatrix}
1 & 0 \nonumber \\
1 & 0 \end{bmatrix}$.
Since $\Prb(a=1) = 1$, the controller doesn't accrue any profit by offering services. Therefore the instantaneous and continuation payoff is $0$.
\begin{equation} \label{eq:P3}
(\Rightarrow) Q(\pi,u) = 0.
\end{equation}
\end{compactenum}
The result follows from (\ref{eq:P1}), (\ref{eq:P2}) and (\ref{eq:P3}). \qed

\textit{\underline{Proof of Theorem~\ref{thm:OPRM}}}:\\
Define the following:
\begin{align*}
{\hspace{-0.2cm}}\delta^* &= \min \{ \pi(2) | ~\eta^{y=1}(2) \geq 1-\alpha\}, \\
{\hspace{-0.2cm}}\gamma^* &= \{ \pi | (u^{H}(\pi)-\beta) \times \textbf{1}^\prime B_{y=2} \pi   + \rho \E V(\pi) = 0\}, \\
{\hspace{-0.2cm}}\pi^*(2) &= \max \{\delta^*,\gamma^* \}, \\
{\hspace{-0.2cm}}\pi^{**}(2) &= \{ \pi(2) | (u^{L}(\pi)-\beta) + \rho V(\pi) = \\  
& (u^{H}(\pi)-\beta) \times \textbf{1}^\prime B_{y=2} \pi   + \rho \E V(\pi)\}.
\end{align*}
\begin{compactenum}
\item[i.)] Consider $\pi(2) \in [0,\pi^*(2))$. We will show that $\{\max~Q(\pi,u) = 0 \}$.\\ 
Let $V(0)$ denote the value at $\pi = \begin{bmatrix} 1 \\ 0 \end{bmatrix}$. As $\pi(2)\rightarrow 0$, we have $\E V(0) \rightarrow V(0)$ and $u^{H}(0) = u^L(0) \rightarrow 1 - \frac{1}{\alpha}$. 

Assume on the contrary $V(\pi) = (u^{L}(\pi)-\beta) + \rho V(\pi)$. As $\pi(2)\rightarrow 0$, $V(0) = \frac{(1 - \frac{1}{\alpha} - \beta)}{(1-\rho)}$. Since $\alpha \in (0,1]$, $\frac{1}{\alpha} \geq 1$ and $V(0)<0$. From Theorem~\ref{thm:RM}, $V(\pi) \geq 0$. Contradiction. \\
Similarly if $V(\pi) = (u^{H}(\pi)-\beta) \times \textbf{1}^\prime B_{y=2} \pi   + \rho \E V(\pi)$, we have $V(0)<0$. Therefore, $V(0)=0$ and we have $Q(0,u^H(0)) < 0$ and $Q(0,u^L(0))<0$. 

From the convexity of the value function, $\E V(\pi) \geq V(\pi)$. Since $Q(\pi,u^H(\pi))<0$ for $\pi(2) =[0,\pi^*(2))$, by definition of $\pi^*(2)$, we have
\begin{align*}
{\hspace{-0cm}} (u^{H}(\pi)-\beta) \times \textbf{1}^\prime B_{y=2} \pi   &+ \rho \E V(\pi) < 0 \\
{\hspace{-0cm}}V(\pi) \geq 0 \rightarrow \E V(\pi) \geq 0 &~~\text{by Jensen's Inequality.} \\
{\hspace{-0cm}}\therefore (u^{H}(\pi)-\beta) &< 0. \\
{\hspace{-0cm}}(u^{H}(\pi)-\beta) < 0 \rightarrow (u^{L}(\pi)-\beta)&<0 ~~\text{from Lemma~\ref{lem:IncPB}}.
\end{align*}
If on the contrary $V(\pi) = (u^{L}(\pi)-\beta) + \rho V(\pi)$, then $V(\pi) < 0$; a contradiction. 
\begin{align*}
\therefore Q(\pi,u^L(\pi)) &<0 ~~\text{for all}~\pi(2) =[0,\pi^*(2)].\\
\Rightarrow V(\pi) &= 0~~\text{for all}~\pi(2) =[0,\pi^*(2)].
\end{align*}

\item[ii.)] $\pi^{**}(2) = \{ \pi(2) | Q(\pi,u^H(\pi) = Q(\pi,u^L(\pi))) \}$. We will show that for $\pi(2) \in (\pi^{**}(2), 1],~Q(\pi,u^L(\pi))) > Q(\pi,u^H(\pi))) > 0$. 

Assume $Q(\pi,u^{H}(\pi))~>Q(\pi,u^{L}(\pi))$ on the contrary. Consider $\pi(2) \rightarrow 1$. Let $V(1)$ and $b( =\textbf{1}^\prime B_{y=2}\pi) \in [0,1]$ denote the values at $\pi =~\begin{bmatrix} 0 \\ 1 \end{bmatrix}$.  We have
\begin{align*}
u^{H}(1) = u^{L}(1) \rightarrow 1 ~&{\text{and}}~\E V(1) \rightarrow V(1)\\
\Rightarrow (1-\beta) \times b  + \rho \E V(1) &> (1-\beta) + \rho V(1) \\
\Rightarrow b &> 1, ~\text{a contradiction as $\beta >0$.} \\
\Rightarrow Q(\pi,u^{L}(\pi)) &> Q(\pi,u^{H}(\pi)).
\end{align*}

From Theorem~\ref{thm:RM}, $V(\pi) \geq 0$ and therefore, $\E V(\pi) \geq 0$. 
\begin{align*}
\text{For}~\pi(2) \in [\pi^{**}(2), 1],~ (u^{H}(\pi)-\beta) \times \textbf{1}^\prime B_{y=2}\pi &> 0 \\
(\Rightarrow) Q(\pi,u^{H}(\pi)) &> 0.
\end{align*}
\item[iii.)] Since $\pi^{**}(2) = \{ \pi(2) | Q(\pi,u^H(\pi) = Q(\pi,u^L(\pi))) \}$ and $Q(\pi,u^L(\pi))<0$ for all $\pi(2) =[0,\pi^*(2)]$ , from part~$(ii)$ we have 
\begin{equation*}
{\hspace{-0.5cm}}Q(\pi,u^H(\pi))) > Q(\pi,u^L(\pi)))~\text{for all}~\pi(2) \in [\pi^*(2),\pi^{**}(2)).
\end{equation*}
Note that $Q(\pi,u^H(\pi)))>0~\text{for all}~\pi(2) \in [\pi^*(2),\pi^{**}(2))$ by definition of $\pi^{*}(2)$ and the fact that $Q(\pi,u) \uparrow \pi$ (Theorem~\ref{thm:VFRM}). \qed
\end{compactenum} 

\textit{\underline{Proof of Theorem~\ref{thm:MPT}}}:\\
The public belief $\pi_k$ is a martingale when the state is a random variable, i.e, $\E[\pi_{k+1}| \mathcal{F}_k] = \pi_k$; see \cite{Cha04,BOOV08}.\\
It can easily be verified{\footnote{Note that the matrix $B$ is TP2. It can be seen that derivative of $u^H(\pi)$ is strictly decreasing and the derivative of $u^L(\pi)$ is strictly increasing with respect to $\pi(2)$ for any $\alpha \in (0,1]$.}} that $u^{H}(\pi)$ is a concave function and $u^{L}(\pi)$ is a convex function of $\pi$ for $\alpha \in (0,1].$
\begin{compactenum}
\item[i.)] \textit{Self-Interested}: For $\epsilon \rightarrow 0$, we have for $\pi_{k}(2),~\pi_{k+1}(2) \in [0,\pi^{**}(2))$, $u_k = u^{H}(\pi_k)$ and it satisfies $\E[u^{H}(\pi_{k+1})| \mathcal{F}_k] \leq u_k$ by Jensen's inequality. 

 We know that $u^{L}(\pi) \leq u^{H}(\pi)$ from Lemma~\ref{lem:IncPB}. For the case of $\pi_{k}(2) \in [\pi^*(2),\pi^{**}(2))$ and  $\pi_{k+1}(2) \in [\pi^{**}(2),1]$, we have 
\begin{equation*}
\hspace{-0.7cm} \E[u_{k+1}| \mathcal{F}_k] = \E[u^{L}(\pi_{k+1})|\mathcal{F}_k] \leq  \E[u^{H}(\pi_{k+1})|\mathcal{F}_k] \leq u_{k}.
\end{equation*}
   
Note that the belief is frozen in $[\pi^{**}(2),1]$, so $\pi_{k+1}(2) \in [\pi^*(2),\pi^{**}(2))$ and  $\pi_{k}(2) \in [\pi^{**}(2),1]$ is irrelevant. 
\item[ii.)] \textit{Altruistic}: Here $\pi^{**}(2) = 1$. For $\epsilon \rightarrow 0$, we have for $\pi_{k}(2),~\pi_{k+1}(2) \in [0,1]$, $u_k = u^{H}(\pi_k)$ and it satisfies $\E[u^{H}(\pi_{k+1})| \mathcal{F}_k] \leq u_k$ by Jensen's inequality. 
\end{compactenum} 


\bibliographystyle{ieeetr}
\bibliography{references}

\end{document}